\documentclass[fleqn]{mat01}
\usepackage{times,mathtimy,amssymb,latexsym}
\begin{document}

\setcounter{page}{411} \firstpage{411}




\renewcommand\theequation{\thesection\arabic{equation}}

\newtheorem{theore}{Theorem}
\renewcommand\thetheore{\arabic{theore}}
\newtheorem{theor}[theore]{\bf Theorem}
\newtheorem{lem}[theore]{\it Lemma}
\newtheorem{pot}[theore]{\it Proof of Theorem}
\newtheorem{coro}[theore]{\rm COROLLARY}
\newtheorem{exam}{\it Example}
\newtheorem{case}{\it Case}

\def\remar{\trivlist \item[\hskip \labelsep{\it Remark.}]}
\def\d{\mbox{\rm d}}

\title{A generalization of d'Alembert formula}

\markboth{Yu-Hsien Chang and Cheng-Hong Hong}{A generalization of
d'Alembert formula}

\author{YU-HSIEN CHANG and CHENG-HONG HONG}

\address{Department of Mathematics, National Taiwan Normal University,
88 Sec.~4, Ting Chou Road, Taipei, Taiwan, Republic of China\\
\noindent E-mail: changyh@math.ntnu.edu.tw; hong838@yahoo.com.tw}

\volume{117}

\mon{August}

\parts{3}

\pubyear{2007}

\Date{MS received 8 June 2006; revised 18 August 2006}

\begin{abstract}
In this paper we find a closed form of the solution for the factored inhomogeneous
linear equation
\begin{equation*}
\prod_{j=1}^{n}\left( \frac{\hbox{d}}{\hbox{d}t}-A_{j}\right) u(t) =f(t).
\end{equation*}
Under the hypothesis $A_{1},A_{2}, \dots, A_{n}$ are infinitesimal generators of
mutually commuting strongly continuous semigroups of bounded linear operators on a
Banach space $X$. Here we do not assume that $A_{j}$s are distinct and we offer the
computational method to get explicit solutions of certain partial differential
equations.\end{abstract}


\maketitle

\section{Introduction}

Many homogeneous differential equations can be transformed to factored equations of
the form:
\begin{equation}
\begin{cases}
\prod_{j=1}^{n}\left(\dfrac{\hbox{d}}{\hbox{d}t}-A_{j}\right) u(t) = 0, \\[.6pc]
u^{(k)}(0) = x_{k}\in X,k=0,1,2, \dots, n-1,
\end{cases}
\end{equation}
where $A_{1},A_{2}, \dots, A_{n}$ are infinitesimal generators of mutually commuting
strongly continuous semigroups of bounded linear operators on a Banach space $X$ and
$u^{(k)} (0)$ is the $k$th derivative of $u(t)$ at $t=0$. Under the hypotheses,

\begin{enumerate}
\renewcommand\labelenumi{(A0)}
\leftskip .31cm
\item $A_{1}, A_{2}, \dots, A_{n}$ generate $C_{0}$-semigroups on a
Banach space $X$ which are mutually commuting. That is
$\hbox{e}^{tA_{j}}\hbox{e}^{sA_{k}}=\hbox{e}^{sA_{k}}\hbox{e}^{tA_{j}}$ for all $t, s
\geq 0$ and for all $j, k\in \{ 1, 2, \dots, n\}$.

\renewcommand\labelenumi{(A1)}
\item $A_{j}-A_{k}$ is injective if $j\neq k$;

\renewcommand\labelenumi{(A2)}
\item range $(A_{j}-A_{k})$ is large enough for $j\neq k$.
\end{enumerate}

Goldstein {\it et~al} \cite{2} found the solution of (1.1) by the d'Alembert formula,
which has the form
\begin{equation}
u (t) =\sum_{j=1}^{n} \hbox{e}^{tA_{j}}x_{j}.
\end{equation}

One may say $u (t)$ is either a strong or mild solution of (1.1) and we ignore this
issue for the moment. Recently, the abstract d'Alembert formula has been broadly
applied to equipartition of energy and scattering theory (see for e.g. \cite{3,4,5});
it also has been extended to semigroups that are not strongly continuous (see
\cite{6}). However in these papers, the authors always assume that all $A_{j}$s are
different. They directly showed that the function $u (t) $ given by (1.2) is a
solution of (1.1) by putting $u (t)$ into the differential equation in their papers.
This motivated us to consider the case that some of the $A_{j}$s in the abstract
factored linear equation~(1.1) are equal. We get a unique solution of (1.1) by a
constructive way. The most interesting point is that one may easily follow this
process to get the explicit form of certain differential equation (see \S4).

Throughout this paper we always assume that hypothesis (A0) holds. Under this
assumption $A_{1}, A_{2}, \dots, A_{n}$ are mutually commuting, therefore permuting
the orders of the operators in (1.1) will not influence the solution. Thus, one may
permute operators in (1.1) such that the same operators put together as

\begin{enumerate}
\renewcommand\labelenumi{(P1)}
\leftskip .31cm
\item $\left.\right.$\vspace{-1.85pc}
\begin{align*}
\hskip -1.25pc A_{1}&=A_{2}=\dots =A_{S_{1}}=B_{1}, A_{S_{1}+1}= A_{S_{1}+2}=\cdots
=A_{S_{1}+S_{2}}=B_{2}, \ldots .\\[.2pc]
\hskip -1.25pc &\quad\,
A_{\big(\sum_{j=1}^{i-1}S_{j}\big)+1}=\cdots
=A_{\sum_{j=1}^{i}S_{j}}=B_{i} \quad \hbox{and} \quad
\sum_{j=1}^{i}S_{j}=n.
\end{align*}
\end{enumerate}
With these notations we assume that

\begin{enumerate}
\renewcommand\labelenumi{(Al)$^{\prime}$}
\leftskip .35cm
\item $B_{j}-B_{k}$ is injective if $j\neq k$;

\renewcommand\labelenumi{(A2)$^{\prime}$}
\item range $(B_{j}-B_{k})$ is large enough for $j\neq k$.
\end{enumerate}

Furthermore, to the inhomogeneous initial value problem
\begin{align}
\prod_{j=1}^{n}\left( \frac{\hbox{d}}{\hbox{d}t}-A_{j}\right) u(t)
&=f(t),\quad f(0) \neq
0,\nonumber\\[.2pc]
u^{(k)}(0) &=x_{k}, k = 0, 1, 2, \dots, n-1,
\end{align}
we also assume that

\begin{enumerate}
\renewcommand\labelenumi{(H1)}
\leftskip .31cm
\item $f\in C^{1}([0,T];X)\cap (\cup _{i=1}^{n}C[0,T];[D(A_{i})])$, where $[D(A_{i})]$
is the Banach space $D(A_{i})$ equipped with the graph norm.
\end{enumerate}
Under these assumptions we have the following results.

\begin{theor}[\!]
Suppose the assumptions ${\rm (A0)}, {\rm (A1)}^{\prime}$ and
${\rm (A2)}^{\prime}$ are all fulfilled$,$ then there exists an
unique solution of the homogeneous initial value problem~$(1.1)$
which can be expressed as
\begin{align}
u(t) &= \sum_{k_{1}=0}^{S_{1}-1}\frac{t^{k_{1}}}{k_{1}!}T_{B_{1}}(t) y(n,k_{1})
+\sum_{k_{2}=0}^{S_{2}-1}\frac{t^{k_{2}}}{k_{2}!} T_{B_{2}}(t) y(n,S_{1}+k_{2})
\nonumber\\[.5pc]
&\quad\,+\cdots+\sum_{k_{i}=0}^{S_{i}-1}\frac{t^{k_{i}}}{k_{i}!}T_{B_{i}}(t)
y\left(n, \left(\sum_{j=1}^{i-1}S_{j}\right)+k_{i}\right),
\end{align}
\end{theor}
where $S_{j}$ is the multiplicity of $B_{j}$, $\sum_{j=1}^{i}S_{j}=n,
\{T_{B_{j}}(t)\}_{t\geq 0}$ is the $C_{0}$-semigroup generated by\ $B_{j}$ and the
relation between the coefficient vector
\begin{equation*}
\vec{y}=(y(n,0), y (n,1), \dots, y(n, n-1))^{T},
\end{equation*}
and the initial data vector
\begin{equation*}
\vec{x}=(u(0), u^{\prime} (0), \dots, u^{(n-1)} (0))^{T}
\end{equation*}
can be represented as $M_{n}^{-1}\vec{x}=\vec{y}$. Here the matrix
$M_{n}$ is composed by the sub-matrices $(B_{j})_{S_{j}}, j = 0,
1, 2, \dots, i$, that is
\begin{equation}
M_{n}=[(B_{1})_{S_{1}} \ (B_{2})_{S_{2}} \ \cdots \ (B_{i})_{S_{i}}]
\end{equation}
and the sub-matrix $(B_{j})_{S_{j}}$ is a $n\times S_{j}$ matrix which is formed by
the first $S_{j}$ columns of the matrix
\begin{equation}
(B_{j}) _{n}=\begin{bmatrix}
I & 0 &  &  &  \\[.1pc]
B_{j} & I & 0 &  &  \\[.1pc]
B_{j}^{2} & 2B_{j} & I &  &  \\[.1pc]
\vdots & \vdots & \vdots & \ddots & 0 \\[.1pc]
B_{j}^{n-1} & C_{n-2}^{n-1}B_{j}^{n-2} & C_{n-3}^{n-1}B_{j}^{n-3} & \cdots &I
\end{bmatrix}.
\end{equation}

\begin{theor}[\!]
Suppose the assumptions ${\rm (A1)}^{\prime}, {\rm (A2)}^{\prime}$
and ${\rm (H1)}$ are fulfilled$,$ $u^{(k)}( 0) =0$ $(k = 0, 1, 2,
\dots,n-1)$ and $f(0) \neq 0$, the nontrivial solution of the
inhomogeneous initial value problem~$(1.4)$ can be expressed as
\begin{align}
u(t) &=\sum_{k_{1}=0}^{S_{1}-1}\int_{0}^{t}\frac{1}{k_{1}!}(t-s)
^{k_{1}}T_{B_{1}}(t-s) Z(n,k_{1}) f(s) \ {\rm d}s\nonumber\\[.3pc]
&\quad\,+\sum_{k_{2}=0}^{S_{2}-1}\int_{0}^{t}\frac{1}{k_{2}!}(t-s)
^{k_{2}}T_{B_{2}}(t-s) Z(n,S_{1}+k_{2}) f(s) \ {\rm d}s+\cdots\nonumber\\[.3pc]
&\quad\,+\sum_{k_{i}=0}^{S_{i}-1}\int_{0}^{t}\frac{1}{k_{i}!}(t-s)
^{k_{i}}T_{B_{i}}(t-s) Z \left(n,\sum_{j=1}^{i-1}S_{j}+k_{i}\right) f(s) \ {\rm
d}s,\end{align}
\end{theor}
where $S_{i}$ is the multiplicity of $B_{i}$ and $\sum S_{i}=n$,
$\{T_{B_{i}}(t)\}_{t\geq 0}$ is the $C_{0}$-semigroup generated by
$B_{i}$. Furthermore, if\ \ $\vec{h}=(0, 0, \dots, I)$, $I$ is the
identity operator on $C^{1}([0,T];X)\cap (\cup
_{i=1}^{n}C[0,T];[D(A_{i}) ]) $, $\vec{z} =(Z(n,0), Z(n,1), \dots,
Z(n,n-1))^{T}$ is a vector and $M_{n}=[(B_{1})
_{S_{1}}(B_{2})_{S_{2}}\cdots (B_{i}) _{S_{i}}]$ is the matrix
defined as in Theorem~1 and the following relation holds:
\begin{equation}
M_{n}^{-1}\vec{h}=\vec{z}.
\end{equation}

Furthermore, if any one of the initial data of (1.3) does not
equal to zero, then the nontrivial solution of the inhomogeneous
initial value problem~(1.3) can be obtained by combining the
results of Theorems~1 and 2 (see Corollary~3).

\section{Homogeneous equation}

For proving Theorem~1, we will use the following three lemmas as
preliminaries. The proof of them are either straightforward or can
be found in ref.~[8], which we omit here.

\setcounter{theore}{0}
\begin{lem}
Let $u_{j+1}(t) =\prod_{k=1}^{j}\big(\frac{{\rm d}}{{\rm d}t}-A_{k}\big) u_{1}(t) $,
$j=1,2, \dots, n-1$ for all $t\geq 0$ and assume that $u_{j}(t) \in $ $D(A_{j})$ for
all $t\geq 0$. Then $(1.1)$ is equivalent to the vector-valued initial value problem
\setcounter{equation}{0}
\begin{equation}
\left\{ \begin{array}{l@{\,}l} \ \dfrac{{\rm d}\vec{u}(t) }{{\rm
d}t}&=\begin{bmatrix}
A_{1} & 1 &  &  &  \\[.1pc]
& A_{2} & 1 & 0 &  \\[.1pc]
&  & \ddots &  &  \\[.1pc]
& 0 &  & A_{n-1} & 1 \\[.1pc]
&  &  &  & A_{n}
\end{bmatrix} \vec{u}(t), \\[3pc]
\vec{u}(0) &=(u_{1}(0), u_{2}(0), \dots, u_{n}(0)).
\end{array}\right.\end{equation}
\end{lem}
where $\vec{u}(t) \equiv (u_{1}(t), u_{2}(t), \dots, u_{n}(t))$.
The initial data component element in (2.1) is
\begin{equation}
\hskip -4pc u_{m}(0) =u_{m}^{\ast }=x_{m-1}+\sum_{k=0}^{m-2}(-1) ^{k}\sum
(A_{i_{j_{1}(k) }}\cdots A_{i_{j_{k}(k) }}x_{m-k}), \quad m = 1, 2, \dots, n
\end{equation}
with $i_{j_{s}(k)}< i_{j_{t}(k)}$ if $s<t$ for all $i_{j_{s}(k)}\in \{1, 2, \dots,
n\}$ and $x_{m}$ is the initial data in (1.1).

\begin{lem}
As long as $x\in D(A_{i}) \cap D(A_{j})$,
\begin{equation*}
\int_{0}^{t}T_{i}(t-s) T_{j}(s)x{\rm d}s=(A_{j}-A_{i})^{-1}(T_{i}(t) -T_{j}(t))x{\rm
d}s
\end{equation*}
for all $1\leq i,j\leq n$ and for all $0<s\leq t$. Furthermore$,$ for any integer
$k\geq 2$,
\begin{align*}
\hskip -4pc \int_{0}^{t}T_{i}(t-s) \frac{s^{k}}{k!} T_{j}(s)x{\rm
d}s &=\frac{t^{k}}{k!}(A_{j}-A_{i}) ^{-1}T_{j}(t) x\\[.4pc]
&\quad\, -\frac{1}{(k-1)!}(A_{j}-A_{i}) ^{-1}\int_{0}^{t}T_{i}(t-s) s^{k-1}T_{j}(s)
x{\rm d}s
\end{align*}
for all $1\leq i,j\leq n$ and for all $0<s\leq t.$
\end{lem}

\begin{lem}
Let $X$ be a reflexive Banach space and let $A$ be the infinitesimal generator of a
$C_{0}$-semigroup $\{T(t) \}_{t\geq 0}$ on $X$. If $f$\ is a Lipschitz continuous
function on $[0,T],$ then for every $x\in D(A) $ the initial value problem
\begin{equation}
\begin{cases}
\dfrac{{\rm d}}{{\rm d}t}u(t)\hskip -.6pc &= Au(t) + f(t) \\[.4pc]
u(0)\hskip -.6pc &= x \end{cases}
\end{equation}
has a unique solution on $\lbrack 0,T\rbrack $ given by
\begin{equation*}
u(t) =T(t) x+\int_{0}^{t}T(t-s) f(s) {\rm d}s.
\end{equation*}
{\rm (}see p.~$109$ of {\rm [8])}.
\end{lem}

\setcounter{theore}{0}
\begin{pot}
{\rm To show that this theorem is true, we prove two special cases at first.}
\end{pot}

\begin{case} {\rm Suppose
\begin{equation}
\prod_{j=1}^{n}\left(\dfrac{{\rm d}}{{\rm d}t}-A_{j}\right) u(t) = \left(\dfrac{{\rm
d}}{{\rm d}t}-A \right) ^{n}u\ (t) = 0
\end{equation}
(that is $A_{j} = A = B_{1}$ for all $j = 1, 2, \dots, n$). By Lemma~1, solving~(1.1)
is equivalent to solving~(2.1). One may find the solution of~(2.1) by successively
solving $u_{k}(t)$, $k = n, n-1, \dots, 1$. In fact $u_{n}$ is the solution of the
initial value problem
\begin{equation*}
\begin{cases}
\left(\dfrac{{\rm d}}{{\rm d}t}-A \right) \ u_{n}(t) = 0\\[.8pc]
u_{n}(0) =\ u_{n}^{\ast}
\end{cases},
\end{equation*}
where $u_{n}^{\ast}$ is defined by~(2.2). Then $u_{n}(t) = T(t) u_{n}^{\ast}$ where
$\{ T(t) \}_{t\geq 0}$ is the $C_{0}$-semigroup generated by $A$. One may get
$u_{n-1}(t)$ from $u_{n}(t)$ by solving the follow equation:\vspace{-.15pc}
\begin{equation*}
\begin{cases}
\dfrac{{\rm d}}{{\rm d}t}u_{n-1}(t) = Au_{n-1}(t) +u_{n}(t) \\[.7pc]
u_{n-1}(0) = u_{n-1}^{\ast}
\end{cases},
\end{equation*}
where $u_{n-1}^{\ast}$ is defined by (2.2).

In fact,\vspace{-.3pc}
\begin{align*}
u_{n-1}(t) &= T(t) u_{n-1}^{\ast }+ \int_{0}^{t}T(t-s) u_{n}(s) \ {\rm d}s\\[.3pc]
&=T(t) u_{n-1}^{\ast} + \int_{0}^{t}T(t-s) T(s) u_{n}^{\ast}{\rm d}s\\[.3pc]
&=T(t) u_{n-1}^{\ast} + tT(t) u_{n}^{\ast} = T(t) y(2,0) +tT(t) y(2,1),
\end{align*}
where $y(2,0) =u_{n-1}^{\ast }$ and $y(2,1) =u_{n}^{\ast }.$ In general, if we get
\begin{equation*}
u_{n-k}(t) =\sum_{j=0}^{k}\frac{t^{j}}{j!} T(t) y(k+1,j) \quad \hbox{for}
\quad k = 0, 1, 2, \dots, n-1,
\end{equation*}
then one may get $u_{n-(k+1) }(t) $ by solving the initial value
problem
\begin{align*}
\begin{cases}
\left(\dfrac{{\rm d}}{{\rm d}t}-A\right) u_{n-(k+1) }(t) =\ u_{n-k}(t)
=\sum_{j=0}^{k}\frac{t^{j}}{j!}T(t) y(k+1,j) ; \\[.3pc]
u_{n-(k+1) }(0) =u_{n-(k+1) }^{\ast}.
\end{cases}
\end{align*}
The solution of this initial value problem is
\begin{align*}
u_{n-(k+1) }(t) &= T(t) \ u_{n-(k+1) }(0) + \int_{0}^{t}T(t-s)
\sum_{j=0}^{k}\frac{t^{j}}{j!}T(s) y(k+1,j) \ {\rm d}s\\[.3pc]
&=T(t) \ u_{n-(k+1) }(0) +\sum_{j=0}^{k}\frac{t^{j+1}}{(j+1)!} T(t)
y(k+1,j).\pagebreak
\end{align*}

\noindent For simplifying the notation, we denote
\begin{align}
u_{n-(k+1) }(t) =\sum_{j=0}^{k+1}\frac{t^{j}}{j!} T(t) y(k+2,j),
\end{align}
where $y(k + 2, 0) = u_{n-(k+1)}(0)$ and $y(k+2,j)$ in (2.5) which is equal to
$y(k+1,j-1)$ in (2.4) for all $j = 1, 2, \ldots, k$. Although, the expression of
vector $y(n, k)$ is not very clear till now, we will find the expression $y(n,k) $ in
terms of initial values $u^{(k)}(0)$s at the end of this proof.}
\end{case}

\begin{case}
{\rm Suppose eq.~(1.1) is expressed with some suitable initial data as
\begin{equation}
\left(\dfrac{{\rm d}}{{\rm d}t}-B\right) \left(\dfrac{{\rm d}}{{\rm
d}t}-A\right)^{n-1}u(t) = 0.
\end{equation}
As in Case~1, one can get the solution of the equation
\begin{equation}
\left(\frac{{\rm d}}{{\rm d}t}-A\right)^{n-1}u(t) = 0.
\end{equation}
Denote the solution of (2.7) by $u_{2}(t) =\sum_{k=0}^{n-2} \frac{t^{k}}{k!}$
$T_{A}(t)y(n-1,k)$, where $\{T_{A}(t)\}_{t\geq 0}$ is the $C_{0}$-semigroup generated
by $A$.

Then solving the initial value problem~(2.6) is equivalent to solving the following
initial value problem
\begin{equation}
\begin{cases}
\left(\dfrac{{\rm d}}{{\rm d}t}-B\right) u_{1}(t) &=u_{2}(t) ; \\[.4pc]
 u_{1}(0) &=u_{1}^{\ast}.
\end{cases}
\end{equation}
By Lemma~2,
\begin{align*}
u_{1}(t) &=T_{B}(t)u_{1}^{\ast } + \int_{0}^{t}T_{B}(t-s) u_{2}(s) \ {\rm d}s\\[.3pc]
&=T_{B}(t)u_{1}^{\ast } + \int_{0}^{t}T_{B}(t-s)
\sum_{k=0}^{n-2}\frac{s^{k}}{k!}T_{A}(s)y(n-1,k) \ {\rm d}s\\[.3pc]
&=T_{B}(t)u_{1}^{\ast }+\int_{0}^{t}T_{B}(t-s) T_{A}(s)y(n-1,0) \ {\rm d}s + \cdots\\[.3pc]
&\quad\,+\int_{0}^{t}T_{B}(t-s) \frac{s^{k}}{k!}T_{A}(s)y(n-1,k) +\cdots \\[.3pc]
&\quad\,+\int_{0}^{t}T_{B}(t-s) \frac{s^{(n-2) }}{(n-2)!}T_{A}(s)y(n-1,n-2)\hbox{d}s\\[.3pc]
&\qquad\,\vdots\\
&=T_{B}(t)u_{1}^{\ast }+\sum_{k=1}^{n-2}\frac{t^{k}}{k!}(A-B)^{-1}T_{A} \ y(n-1,k)\\[.4pc]
&\quad\,+\sum_{k=2}^{n-2}(-1) ^{1}\frac{t^{k-1}}{(k-1)!}(A-B)^{-2}T_{A} \ y(n-1,k) +
\cdots
\end{align*}
\begin{align*}
\phantom{u_{1}(t)}&\quad\,+(-1)^{n-3} t(A-B) ^{-(n-2)}T_{A} y(n-1,n-2)\\[.3pc]
&\quad\, +\sum_{k=0}^{n-1}(-1) ^{k} (A-B)^{-(k+1)}(T_{A}-T_{B})y(n-1,k),
\end{align*}
where $\{T_{B}(t)\}_{t\geq 0}$\ is \ $C_{0}$-semigroup generated by $B$.

Rewrite $u_{1}$ in terms of increasing degree of $t$, one may get
\begin{align}
\begin{split}
\hskip -1.25pc u_{1}(t) &=T_{B}(t)\left[u_{1}(0) +\sum_{k=0}^{n-2}(-1) ^{k+1}(A-B)
^{-(k+1)} y(n-1,k)
\right]\\[.4pc]
&\quad\,+ T_{A}(t) \left[\sum_{k=0}^{n-2}(-1) ^{k}(A-B) ^{-(k+1)} y(n-1,k) \right]\\[.4pc]
&\quad\,+tT_{A}(t) \left[\sum_{k=1}^{n-2}(-1) ^{k-1}(A-B) ^{-k} y(n-1,k) \right]
+\cdots\\[.4pc]
&\quad\,+\frac{t^{j}}{j!}T_{A}(t) \left[\sum_{k=j}^{n-2}(-1) ^{k-j}(A-B) ^{-(k+1-j) }
y(n-1,k) \right] +\cdots\\[.4pc]
&\quad\,+\frac{t^{n-1}}{(n-1)!}T_{A}(t) [ (-1) ^{0}(A-B) ^{-1}y(n-1,n-2)].
\end{split}
\end{align}
For simplifying the notation, we denote
\begin{align}
\begin{split}
u_{1}(t) &=T_{B}(t)[y(n,0)]+T_{A}(t)[y(n,1) ]+tT_{A}(t)[y(n,2) ]+\cdots\\[.3pc]
&\quad\,+\frac{t^{k}}{k!}T_{A}(t)[y(n,k+1) ]+\cdots +\frac{t^{n-1}}{(n-1) !}T_{A}(t)
[y(n,n-1)], \end{split}
\end{align}
where \begin{align*}
[y(n,0)] &=\left[ u_{1}^{\ast }+\sum_{k=0}^{n-2}(-1) ^{k+1}(A-B)
^{-(k+1)}y(n-1,k) \right],\\[.4pc]
[y(n,1)] &=\left[ \sum_{k=0}^{n-2} (-1) ^{k}(A-B) ^{-(k+1)}y(n-1,k) \right],\\[.4pc]
[y(n,j)] &=\left[ \sum_{k=j}^{n-2} (-1) ^{k-j}(A-B) ^{-(k+1-j) }y(n-1,k) \right]\\[.4pc]
&\qquad (2\leq j\leq n-2) \quad \hbox{and}\\[.4pc]
[y(n,n-1)] &= [(-1)^{2(n-1) }(A-B)^{-1}y(n-1,n-2)].
\end{align*}}
\end{case}

For the general case, one may first permute operators in (1.1)
such that the same operators are put together and then apply the
results in Cases~1 and 2 alternately. Finally, one can reach the
conclusion that the solution of (1.1) can be expressed as in
(1.4). This theorem will be proved as long as one finds the
relation between initial data vector $\vec{x}=(u(0),
u^{\prime}(0), \dots, u^{n-1}(0))^{T}$ and the vector
$\vec{y}=(y(n,0), y(n,1), \dots, y(n,n-1))^{T}$. Since the
solution $u(t)$ can be expressed as the combination of the
terms\vspace{-.3pc}
\begin{align*}
u(t;k,l) = \frac{t^{k}}{k!} \ {\rm e}^{tA_{l}}x, \quad x\in D(A)\vspace{-.3pc}
\end{align*}
we consider the derivatives of $u(t;k,l) .$ From the fact
$C_{m}^{n}=C_{m-1}^{n-1}+C_{m}^{n-1}$, one may get the $i$-th derivative of $u(t;k,l)$
as\vspace{-.3pc}
\begin{align}
\hskip -4pc u(t;k,l)^{(i)} &=\frac{1}{(k-i)!}t^{(k-i)} \ {\rm
e}^{tA_{l}}x+\frac{1}{(k-(i-1) )
!}C_{1}^{i}t^{(k-(i-1)) }A_{l} \ {\rm e}^{tA_{l}} x + \cdots\nonumber\\[.7pc]
\hskip -4pc &\quad\,+\frac{1}{(k-(i-j)) !}C_{j}^{i}t^{(k-(i-j)) }A_{l}^{j} \ {\rm
e}^{tA_{l}}x+\cdots \nonumber\\[.7pc]
&\quad\,+ \frac{1}{k!}C_{i}^{i} \ t^{k}A_{l}^{i} \ {\rm
e}^{tA_{l}}x \quad \hbox{for} \quad i<k
\end{align}
and
\begin{align}
\hskip -4pc u(t;k,l)^{(k+j)} &=C_{j}^{k+j} A_{l}^{j} \ {\rm
e}^{tA_{l}}+\frac{1}{1!}C_{j+1}^{k+j} A_{l}^{j+1}t \ {\rm e}^{tA_{l}} + \cdots + \frac{1}{h!}C_{j+h}^{k+h} A_{l}^{j+h}t^{h} \ {\rm e}^{tA_{l}} +\cdots\nonumber\\[.7pc]
\hskip -4pc &\quad\, +\frac{1}{k!} C_{k+j}^{k+j} A_{l}^{k+j}t^{k}
\ {\rm e}^{tA_{l}} \quad \hbox{for} \quad j=i-k\geq
0.\vspace{-.3pc}
\end{align}
Since we want to find the relation between the coefficient vector $%
\vec{y}$ and initial data vector $\vec{x}$, we need only to
consider the special situation $t=0$. According to (2.11) and
(2.12), if one just considers a single operator $A_{i}$, one can
get\vspace{-.3pc}
\begin{equation}
(A_{i}) _{n}= \begin{bmatrix}
I & 0 &  &  &  \\[.1pc]
A_{i} & I & 0 &  &  \\[.1pc]
A_{i}^{2} & 2A_{i} & I &  &  \\[.1pc]
\vdots & \vdots & \vdots & \ddots & 0 \\[.1pc]
A_{i}^{n-1} & C_{n-2}^{n-1}A_{i}^{n-2} & C_{n-3}^{n-1}A_{i}^{n-3} & \cdots & I
\end{bmatrix}.
\end{equation}
In general, after we permute the operators in (1.1) such that the same operators are
put together as (P1), the equation can be represented as
\begin{align*}
\prod_{i=1}^{S_{1}}\left(\frac{\d}{\d t}%
-B_{1}\right) \prod_{i=1}^{S_{2}}\left(\frac{\d}{\d t}-B_{2}\right) \dots
\prod_{i=1}^{S_{i}}\left( \frac{\d}{\d t}-B_{i}\right) u=0
\end{align*}
and one can find the sub-matrices $M_{n}^{j}$ corresponding to
$\prod_{i=1}^{S_{j}}\big(\frac{{\rm d}}{{\rm d} t}-B_{j}\big) $ as
(1.6). Combining these $i$ sub-matrices together, one may get
$M_{n}$ as (1.5). The uniqueness of the solution follows from
Lemmas~1 and 3 immediately.\pagebreak

\begin{remar}
We use the following example to demonstrate how to get the matrix $M_{n}$ in
Theorem~1. Consider the initial value problem
\begin{align}
\begin{cases}
\left( \dfrac{\d}{\d t}-A\right) \left( \dfrac{\d}{\d t}-A\right) \left( \dfrac{\d}{\d t}%
-B\right) \left( \dfrac{\d}{\d t}-B\right) \left(
\dfrac{\d}{\d t}-C\right) u\left( t\right) =0; \\[.7pc]
u^{(k)}(0) =x_{k}, \quad (k=0,1,2,3,4).%
\end{cases}
\end{align}
\end{remar}

\noindent By comparing with the proof of theorem 1, one may rewrite the differential
equation and initial data in (2.14) as
\begin{align*}
u_{5}(t) =\left( \frac{\d}{\d t}%
-A\right) \left( \frac{\d}{\d t}-B\right) \left( \frac{\d}{\d t}-B\right) \left(
\frac{\d}{\d t}-C\right) u(t)
\end{align*}
and
\begin{align*}
u_{5}^{\ast } =u_{5}(0) &=x_{4}-(A+2B+C)
x_{3}+ (B^{2}+2AB+2BC+AC) x_{2}\\[.1pc]
&\quad\,-(AB^{2}+CB^{2}+2ABC) x_{1}+AB^{2}Cx_{0}.
\end{align*}
At first, solve the initial value problem
\begin{equation*}
\begin{cases}
\left( \dfrac{\d}{\d t}-A\right) u_{5}(t) =0 \\[.6pc]
u_{5}(0) =u_{5}^{\ast }.
\end{cases}
\end{equation*}
One may have
\begin{equation*}
u_{5}(t) =T_{A}(t)u_{5}^{\ast}  =T_{A}(t)y(1,0),
\end{equation*}
where $\{T_{A}(t)\}_{t\geq 0}$ is $C_{0}$-semigroup generated by
$A$. If we rewrite eq.~(2.14) as Lemma~1 and find the solution by
successively solving $u_{k}(t)$ for $k=4,3,2$, then we get
\begin{align*}
u_{2}(t) = T_{B}y(4,0) +tT_{B}y (4,1) +T_{A}y(4,2) + tT_{A}y(4,3).
\end{align*}
The coefficient vector of $u_{2}(t) $ is $\vec{y}_{4}%
=(y(4,0), y(4,1), \dots ,y (4,3)) ^{T}$. By a similar method, we
get the coefficient vector of $u(t) =u_{1}(t)$ as
$\vec{y}_{5}=(y(5,0), y(5,1), \dots,$ $y(5,4))^{T}$ since
$\vec{y}_{4}$ is just a 4-dimensional vector. For vector of
$\vec{y}_{5}$, we need to extend $\vec{y}_{4}$ as a 5-dimensional
vector. We add $u(0)$ into $\vec{y}_{4}$ as the first component of
$\vec{y}_{4}^{\ast}$, that is $\vec{y}_{4}^{\ast}=( u( 0) ,y(4,0),
y(4,1),\dots,y (4,3)) ^{T}$.

Rewrite (2.10) in the matrix form
\begin{align}
\hskip -4.2pc\begin{bmatrix}
y( 5,0) \\[.1pc]
y( 5,1) \\[.1pc]
y( 5,2) \\[.1pc]
y( 5,3) \\[.1pc]
y( 5,4)%
\end{bmatrix}  =\begin{bmatrix}
I & -( B-C) ^{-1} &( B-C) ^{-2} & -(A-C)^{-1} &( A-C) ^{-2} \\[.2pc]
0 &( B-C) ^{-1} & -( B-C) ^{-2} &0 &0 \\[.2pc]
0 &0 &( B-C) ^{-1} &0 &0 \\[.2pc]
0 &0 &0 &( A-C) ^{-1} & -( A-C) ^{-2} \\[.2pc]
0 &0 &0 &0 &( A-C) ^{-1}
\end{bmatrix}
\begin{bmatrix}
u( 0) \\[.1pc]
y( 4,0) \\[.1pc]
y( 4,1) \\[.1pc]
y( 4,2) \\[.1pc]
y( 4,3)%
\end{bmatrix}.
\end{align}
From Theorem~1, the relation between $\vec{y}_{4}^{\ast}=( u( 0)
,y( 4,0) ,y( 4,1),\dots,y( 4,3) ) ^{T}$ and $\vec{u}^{\prime}=( u(
0) ,u_{2}( 0) ,u_{2}^{\prime }( 0) ,u_{2}^{\prime \prime }( 0)
,u_{2}^{( 3) }( 0) ) ^{T}$ can be represented as
\begin{equation}
\begin{bmatrix}
u( 0) \\[.1pc]
y( 4,0) \\[.1pc]
y( 4,1) \\[.1pc]
y( 4,2) \\[.1pc]
y( 4,3)%
\end{bmatrix} =
\begin{bmatrix}
I & 0 & 0 & 0 & 0 \\[.1pc]
0 & I & 0 & I & 0 \\[.1pc]
0 & B & I & A & I \\[.1pc]
0 & B^{2} & 2B & A^{2} & 2A \\[.2pc]
0 & B^{3} & 3B^{2} & A^{3} & 3A^{2}%
\end{bmatrix}^{-1} \
\begin{bmatrix}
u( 0) \\[.2pc]
u_{2}( 0) \\[.2pc]
u_{2}^{\prime }( 0) \\[.3pc]
u_{2}^{\prime \prime }( 0) \\[.3pc]
u_{2}^{( 3) }( 0)%
\end{bmatrix}.
\end{equation}
Followed from Lemma~1, the relation between the initial data of $u( t)$ and $u_{2}(
t)$ can be expressed as
\begin{equation}
\begin{bmatrix}
u( 0) \\[.2pc]
u^{\prime }( 0) \\[.2pc]
u^{\prime \prime }( 0) \\[.2pc]
u^{(3)}( 0) \\[.2pc]
u^{(4)}( 0)%
\end{bmatrix} =
\begin{bmatrix}
I & 0 & 0 & 0 & 0 \\[.1pc]
-C & I & 0 & 0 & 0 \\[.1pc]
0 & -C & I & 0 & 0 \\[.1pc]
0 & 0 & -C & I & 0 \\[.1pc]
0 & 0 & 0 & -C & I%
\end{bmatrix}^{-1} \
\begin{bmatrix}
u( 0) \\[.2pc]
u_{2}( 0) \\[.2pc]
u_{2}^{\prime }( 0) \\[.3pc]
u_{2}^{\prime \prime }( 0) \\[.3pc]
u_{2}^{( 3) }( 0)%
\end{bmatrix}.
\end{equation}
According to (2.15), (2.16) and (2.17), we get
\begin{equation*}
M_{5}=
\begin{bmatrix}
I & I & 0 & I & 0 \\[.1pc]
C & B & I & A & I \\[.1pc]
C^{2} & B^{2} & 2B & A^{2} & 2A \\[.1pc]
C^{3} & B^{3} & 3B^{2} & A^{3} & 3A^{2} \\[.1pc]
C^{4} & B^{4} & 4B^{3} & A^{4} & 4A^{3}%
\end{bmatrix}.
\end{equation*}

\section{Inhomogeneous equation}

In this section, we consider the inhomogeneous initial value problem
\setcounter{equation}{0}
\begin{equation}
\begin{cases}
\prod_{j=1}^{n}\left( \dfrac{\d}{\d t}-A_{j}\right) u( t) =f(
t),\quad f( 0) \neq 0, \\[.8pc]
u^{(k)}( 0) =x_{k},k=0,1,2,\dots ,n-1.%
\end{cases}
\end{equation}
According to superposition principle, one may obtain the solution of (3.1) by
combining the solution of homogeneous initial value problem (1.1) with nonzero initial
data and the solution of inhomogeneous case with zero initial data.

\begin{pot}
{\rm As in the proof of Theorem~1, we consider $A_{j}=A$ for all $j=1,2,\dots,n$
at first. In this case, one may follow the process shown in Theorem~1 to get%
\begin{equation}
\hskip -4pc v_{n-k}( t) =\int_{0}^{t}\frac{1}{k!}( t-s) ^{k}T_{A}( t-s) Z( n,k) f(
s)\, \d s,  \quad (k=0,1,2,\dots, n-1),
\end{equation}
where $v_{k}( t)$ $(k=0,1,2,\dots,n-1)$ is the solution of the initial value problem
\begin{equation*}
\begin{cases}
\left(\dfrac{\d}{\d t}-A\right) v_{n-k}( t) =v_{n-(k-1)}(
t) \\[.8pc]
v_{n-k}( 0) =0%
\end{cases}
\end{equation*}
with $v_{n}( t) =f( t)$. Finally, one may obtain $v_{1}( t) $ to be the solution of
(3.1)
with zero initial data. However, for finding the operators $%
Z( n,k) $s in the representation of (1.7), one may rewrite the initial value problem
(3.1) (with zero initial data) in the equivalent system\break form
\begin{equation}
\begin{cases}
\dfrac{{\rm d}\vec{u}( t) }{{\rm d}t}=
\begin{bmatrix}
A_{1} & 1 &  &  &  \\
& A_{2} & 1 &  &  \\
&  & \ddots & 0 &  \\
& 0 &  & A_{n-1} & 1 \\
&  &  &  & A_{n}%
\end{bmatrix}
\vec{u} (t) +
\begin{bmatrix}
0 \\
0 \\
0 \\
\vdots \\
f%
\end{bmatrix}; \\[1.5cm]
\vec{u}( 0) =\vec{u}_{0}( 0) =( u_{1}( 0) ,u_{2}( 0) ,\dots,u_{n}(
0) ) =( 0,0,\dots,0) ,%
\end{cases}
\end{equation}
where $u_{1}( t) $ is the solution of the initial value problem (3.1) with zero
initial data and
\begin{align*}
u_{j+1}( t) =\prod_{k=1}^{j}\left( \frac{\d}{\d t}-A_{k}\right) u_{1}( t) \in D(
A_{j+1}) \ \hbox{for all} \ t\geq 0,\quad 0\leq j<n.
\end{align*}
We will show that relation (1.8) holds at the end of this theorem.

Further, we consider $A_{1}=B$, $A_{2}=A_{3}=\dots =A_{n}=A$ (i.e. there are only two
distinct operators in eq.~(3.1)).

Following a similar procedure as shown in Theorem~1, we get
\begin{align}
\hskip -4pc u_{n-k}( t) =\int_{0}^{t}\frac{1}{k!}%
( t-s) ^{k}T_{A}( t-s) Z( n,k) f( s)\, \d s,\quad
(k=0,1,2,\dots,n-1).
\end{align}
According to (3.4), the solution of
\begin{equation*}
\prod_{k=2}^{n}\left( \frac{\d}{\d t}-A\right) u( t) =f(t)
\end{equation*}
is
\begin{equation*}
u_{2}( t) =\int_{0}^{t}\frac{1}{( n-2) !}( t-s) ^{n-2}T_{A}( t-s) Z( n,n-2) f( s)\, \d
s.
\end{equation*}
To get the solution of the equation
\begin{equation}
\left( \frac{\d}{\d t}-B\right) \prod_{k=2}^{n}\left( \frac{\d%
}{\d t}-A\right) u( t) =f( t)
\end{equation}
one need only to solve the equation
\begin{equation}
\left(\frac{\d}{\d t}-B\right) u_{1}=u_{2}.
\end{equation}
By Fubini's theorem and integration by parts, it is easy to see that
\begin{align}
u_{1}( t) &=\int_{0}^{t}\frac{1}{( n-1) !}%
( t-\tau ) ^{n-1}( A-B) ^{-1}T_{A}( t-\tau ) Z( n,n-2) f( \tau )\,
\d\tau\nonumber\\[.4pc]
&\quad\,-\int_{0}^{t}\int_{\tau }^{t}\frac{1}{(n-2)!}( s-\tau )
^{n-2}T_{A}( s-\tau )\nonumber\\[.4pc]
&\quad\,\times T_{B}( t-s) ( A-B) ^{-1}Z( n,n-2) f( \tau )\, \d s \d\tau.
\end{align}
Since the operators $A_{j}$s are mutually commuted, without loss of generality we may
assume that the same operator $
\big(\frac{{\rm d}}{{\rm d}t}%
-A_{j}\big)$ in (3.1) are put together such as (P1) and rewrite the differential
equation in (3.1) as
\begin{equation}
\prod_{j=1}^{S_{1}}\left( \frac{\d}{\d t}%
-B_{1}\right)  \prod_{j=1}^{S_{2}}\left( \frac{\d}{\d t}-B_{2}\right) \dots
\prod_{j=1}^{S_{i}}\left( \frac{\d}{\d t}-B_{i}\right) u=f.
\end{equation}
We denote
\begin{equation}
W( t;j,B_{l}) =\sum_{j=0}^{S_{l}-1}\int_{0}^{t}%
\frac{1}{j!}( t-\tau ) ^{j}T_{B_{l}}( t-\tau ) Z( j,l,B_{l}) g( \tau )\, \d \tau
\end{equation}
to be the solution of
\begin{equation}
\prod_{j=1}^{S_{l}}\left( \frac{\d}{\d t}-B_{l}\right) u( t) =g,
\end{equation}
where the index $j$ in $W( t;j,B_{l}) $ denotes the counting number of steps in the
iterative procedure for solving the problem
from the beginning. Then the solution of (3.8) can be represented as%
\begin{align}
u( t) &= \sum_{j=0}^{S_{1}-1} \int_{0}^{t}\frac{1}{j!}( t-\tau )
^{j}T_{B_{1}}( t-\tau ) Z( n,j,B_{1}) f( \tau )\, \d\tau\nonumber\\[.4pc]
&\quad\, + \sum_{j=0}^{S_{2}-1}\int_{0}^{t} \frac{1}{j!}( t-\tau )
^{j}T_{B_{2}}( t-\tau ) Z( n,j,B_{2}) f( \tau )\, \d\tau + \dots\nonumber\\[.4pc]
&\quad\,+\sum_{j=0}^{S_{i}-1}\int_{0}^{t}%
\frac{1}{j!}( t-\tau ) ^{j}T_{B_{l}}( t-\tau ) Z( n,j,B_{i}) f( \tau )\, \d\tau.
\end{align}
Simplifing the notation we denote the coefficient of (3.11) as%
\begin{align*}
Z( n,j,B_{1}) &= Z( n,k_{1})  \ \hbox{for} \ 0\leq j=k_{1}\leq
S_{1}-1;\\[.2pc]
Z( n,j,B_{2}) &= Z( n,S_{1}+k_{2}) \ \hbox{for} \ 0\leq
j=k_{2}\leq S_{2}-1;\\[.2pc]
&\quad\,\vdots\\[.2pc]
Z( n,j,B_{l}) &=Z\left(
n,\sum_{j=1}^{i-1}S_{j}+k_{i}\right) \ \hbox{for} \ 0\leq j=k_{i}\leq S_{i}-1.%
\end{align*}
This theorem will be proved as long as one find the explicit form
of $ Z( n,j,B_{l}) $. One may follow the procedure in Theorem~1 to
find the matrix $M_{n}$. As in Theorem~1, we begin with the
special case that all $A_{i}$s are equal to $A$, then combine the
results of distinct $B_{j}$s part to get the general form $M_{n}$.
When all $A_{i}$s are equal to $A$, the equation is of the form
\begin{equation*}
\prod_{k=1}^{n}\left( \frac{\d}{\d t} -A\right) u=f,
\end{equation*}
and then the solution of this equation is
\begin{equation}
u( t) =W( n,t,A) =\sum_{j=0}^{n-1}\int_{0}^{t}\frac{1}{j!}( t-s) ^{j}T_{A}( t-s) Z(
n,j) f( s)\, \d s.
\end{equation}
One may get (1.8) by continuously differentiating (3.12). $\ $%
The first derivative of (3.12) gives
\begin{align}
u^{\prime }( t) &=\sum_{j=1}^{n-1}\int_{0}^{t} \frac{1}{(j-1)!}(
t-s) ^{(j-1)}T_{A}( t-s) Z( n,j) f( s)\, \d s\nonumber\\[.4pc]
&\quad\,+\sum_{j=0}^{n-1}\int_{0}^{t}\frac{1}{ j!}( t-s) ^{j}AT_{A}( t-s) Z( n,j) f(
s)\, \d s+Z( n,0) f( t).
\end{align}
The initial condition $u^{\prime }( 0) =0$ implies $Z( n,0) f( 0)
=0$. Since $f( 0) \neq 0$, it enforces $Z( n,0) =0$. Continuing
this procedure, one may get
\begin{align*}
u^{( i) }( t) &=\sum_{j=i}^{n-1}\int_{0}^{t}\frac{1}{%
(j-i)!}( t-s) T_{A}( t-s) Z( n,j) f( s)\, \d s\\[.4pc]
&\quad\, +\sum_{j=i-1}^{n-1}\int_{0}^{t}\frac{1}{(j-(i-1))!}%
( t-s) ^{(j-( i-1) )}\\[.4pc]
&\quad\, \times C_{1}^{i}AT_{A}( t-s) Z( n,j) f( s)\, \d
s+\cdots\\[.4pc]
&\quad\, +\sum_{j=i-k}^{n-1}\int_{0}^{t}\frac{1}{(j-( i-k) )!}(
t-s) ^{(j-( i-k) )}\\[.4pc]
&\quad\, \times C_{k}^{i}A^{k}T_{A}( t-s) Z( n,j) f( s)\, \d
s+\cdots\\[.4pc]
&\quad\, +\sum_{j=0}^{n-1}\int_{0}^{t}\frac{1}{j!}( t-s)
^{j}C_{i}^{i}A^{i}T_{A}( t-s)\\[.4pc]
&\quad\, \times  Z( n,j) f( s)\, \d s+Z( n,i-1) f(
t)\\[.4pc]
&\quad\, +C_{1}^{i-1}AZ( n,i-2) f( t)
+\cdots+C_{i-1}^{i-1}A^{i-1}\\[.4pc]
&\quad\, \times Z( n,0) f( t) \quad \hbox{for} \ 1\leq i\leq n
\end{align*}
and
\begin{equation}
\hskip -4pc Z( n,i-1) +C_{1}^{i-1}AZ( n,i-2) +\dots +C_{i-1}^{i-1}A^{i-1}Z( n,0) =0
\quad \hbox{for} \ 1\leq i\leq n-1.
\end{equation}
Finally, one can put $u( t) $, $u^{( i) }( t) $ for $1\leq i\leq n$\ into (3.1) to get
\begin{equation}
\{Z( n,n-1) +C_{1}^{n-1}AZ( n,n-2) +\dots +C_{n-1}^{n-1}A^{n-1}Z( n,0) \} f=f.
\end{equation}
Put (3.14) and (3.15) together and let $\vec{h}=( 0,0,\dots ,I)
^{T}$,\ where$\ I$ is the identity operator on $C^{1}(\lbrack
0,T\rbrack ;X)\cap ( \cup _{i=1}^{n}C\lbrack 0,T\rbrack ;\lbrack
D( A_{i}) \rbrack ) $. Then one can write them in the matrix form
$M_{n}^{-1}\vec{h}=\vec{z}$.\newline

One can apply the superposition principle to get the solution of (3.1) with nonzero
initial data. It is the sum of the solutions of (1.1) and (3.1) with zero initial
data. We summarize this result as follows.}
\end{pot}

\begin{coro}\label{cor1}$\left.\right.$\vspace{.5pc}

\noindent Under the hypotheses of Theorems~$1$ and $2,$ eq.~$(3.1)$ with nonzero
initial data has a solution $u( t) $ which is represented as
\begin{align*}
u( t) &=\sum_{k_{1}=0}^{S_{1}-1}\frac{t^{k_{1}}}{k_{1}!} T_{A_{1}}( t) y( n,k_{1})
+\sum_{k_{2}=0}^{S_{2}-1}
\frac{t^{k_{2}}}{k_{2}!}T_{B_{2}}( t) y( n,S_{1}+k_{2}) +\cdots\\[.4pc]
&\quad\, +\sum_{k_{i}=0}^{S_{i}-1}\frac{t^{k_{i}}}{k_{i}!}%
T_{B_{i}}( t) y\left( n,\left(\sum_{j=1}^{i-1}S_{j}\right)+k_{i}\right) +\cdots \\[.4pc]
&\quad\, +\sum_{k_{1}=0}^{S_{1}-1}\int_{0}^{t}\frac{1}{k_{1}!}(
t-s) ^{k_{1}}T_{B_{1}}( t-s) Z( n,k_{1}) f( s)\, \d s\\[.4pc]
&\quad\, +\sum_{k_{2}=0}^{S_{2}-1}\int_{0}^{t}\frac{1}{k_{2}!}(t-s)
^{k_{2}}T_{B_{2}}( t-s) Z( n,S_{1}+k_{2}) f( s)\, \d s+\cdots\\[.4pc]
&\quad\, +\sum_{k_{i}=0}^{S_{i}-1}\int_{0}^{t}\frac{1}{k_{i}!}( t-s)
^{k_{i}}T_{B_{i}}( t-s) Z\left( n,\sum_{j=1}^{i-1}S_{j}+k_{i}\right) f( s)\, \d s.
\end{align*}
\end{coro}

\section{Applications}

\begin{exam}
{\rm We consider the following initial value problem \setcounter{equation}{0}
\begin{equation}
\begin{cases}
u_{tt}( t,x) +a_{1}u_{tx}( t,x) +a_{2}u_{xx}( t,x) =f( t,x) ,\quad
( t,x) \in [ 0,T]
\times R; \\
u( 0,x) =\phi _{1}( x) ,\quad x\in R\ ; \\
u_{t}( 0,x) =\phi _{2}( x) ,\quad x\in R;%
\end{cases}
\end{equation}
where $a_{1},$\ $a_{2}$ are given constants. Let $E=L^{2}(
R) $, $D( A^{k}) =W_{2}^{k}( R) $ \ and \ $%
A^{k}f=\frac{{\rm d}^{k}}{{\rm d}x^{k}}f$, $k=1,2$ for every $f\in
D( A^{k}) $. Under these notations, (4.1) is equivalent to the
following initial value problem:
\begin{equation}
\begin{cases}
U^{\prime \prime }( t) +a_{1}AU^{\prime }( t) +a_{2}A^{2}U( t) =F(
t), \quad t\in [ 0,T] ;\\[.2pc]
U( 0) =\phi _{1}( x),\quad x\in R; \\[.2pc]
U^{\prime }(0) =\phi _{2}( x),\quad x\in R;%
\end{cases}
\end{equation}
where $U( t) \in L^{2}( R) $ and $F( t) \in L^{2}( R) $ satisfy $(
U( t) ) x=u( t,x) $, $( F( t) ) x=f( t,x) $ for all $( t,x) \in [
0,T] \times R$. It is well-known that $A$\ generates a
$C_{0}$-semigroup $\{ T(
t) \} _{t\geq 0}$ on the Banach space $E$ which satisfies $%
( T( t) f) ( x) =f( x+t) $ for all $f\in L^{2}( R) $ and for all
$( t,x) \in [ 0,T] \times R$ (see for e.g. Ch.~22, item 22.5 of
\cite{7}). Moreover, $T( t) $ can be extended to a $C_{0}$-group.
It is easy to see that
\begin{align*}
\frac{\d^{k}}{\d t^{k}}A^{2-k}h( t) =A^{2-k} \frac{\d^{k}}{\d
t^{k}}h( t),\quad k=1,2 \ \hbox{for any} \ h\in (C^{2} [ 0,T]
:D(A^{2}) ).
\end{align*}
If the characteristic equation $P( z) =z^{2}+a_{1}z+a_{2}=0$ of
(4.2) has a root $z_{1}$ with multiplicity 2, then $z_{1}A$ generates a $%
C_{0}$-semigroup $T_{1}( t) $ which satisfies $( T_{1}( t) f) ( x) =f( x+z_{1}t) $ for
all $f\in L^{2}( R) $ and for all $( t,x) \in [ 0,T] \times R$.

If $\phi _{1}\in D( A^{3}) $ \ and $\phi _{2}\in D( A^{2}) $, then (4.2) can be
rewritten as
\begin{equation}
\begin{cases}
\left(\dfrac{\d}{\d t}-z_{1}A\right) ^{2}U=F( t);\\[.6pc]
U( 0) =\phi _{1}( x);\\[.2pc]
U^{\prime }( 0) =\phi _{2}( x).
\end{cases}
\end{equation}
By Corollary~3, (4.3) has a solution of the form
\begin{align*}
U( t) &=T_{1}( t) y( 2,0) +tT_{1}( t) y( 2,1) +\int_{0}^{t}T_{1}(
t-s) Z( 2,0) F( s)\, \d s\\[.3pc]
&\quad\, +\int_{0}^{t}( t-s) T_{1}( t-s) Z( 2,1) F( s)\, \d s,
\end{align*}
where $y( 2,0) ,$ $y( 2,1) $ and $Z( 2,0) ,$ $Z( 2,1) $ satisfy
\begin{equation*}
\begin{bmatrix}
1 & 0 \\[.1pc]
z_{1}A & 1%
\end{bmatrix}
\begin{bmatrix}
y( 2,0) \\[.1pc]
y( 2,1)%
\end{bmatrix}=
\begin{bmatrix}
\phi _{1} \\[.1pc]
\phi _{2}%
\end{bmatrix};\quad
\begin{bmatrix}
1 & 0 \\[.1pc]
z_{1}A & 1%
\end{bmatrix}
\begin{bmatrix}
Z( 2,0) \\[.1pc]
Z( 2,1)%
\end{bmatrix} =
\begin{bmatrix}
0 \\
I%
\end{bmatrix}.
\end{equation*}
This implies that
\begin{equation*}
U( t) =T_{1}( t) \phi _{1}+tT_{1}( t) ( \phi _{2}-z_{1}A\phi _{1}) +\int_{0}^{t}( t-s)
T_{1}( t-s) F( s)\, \d s.
\end{equation*}
Thus, (4.1) has a solution
\begin{align*}
u( t,x) &=\phi _{1}( x+z_{1}t) +t\phi _{2}( x+z_{1}t) -z_{1}t\phi
_{1}^{\prime }( x+z_{1}t)\\[.4pc]
&\quad\,+\int_{0}^{t}( t-s) f( s,z_{1}( t-s) +x)\, \d s.
\end{align*}}
\end{exam}

\begin{exam}
{\rm We consider the following initial boundary value problem\newline
\begin{equation}
\hskip -4pc\begin{cases}
\dfrac{\partial ^{2}}{\partial t^{2}}( u( t,x) ) +b_{1}%
\dfrac{\partial }{\partial t}( \Delta u( t,x) ) +b_{2}\Delta
^{2}u( t,x) =f( t,x),\quad (
t,x) \in (0,T)\times \Omega; \\[.6pc]
u( 0,x) =\Psi _{1}( x),\quad x\in \Omega; \\[.6pc]
\dfrac{\partial }{\partial t}( u( 0,x) ) =\Psi _{2}( x),\quad
x\in \Omega; \\[.6pc]
u( t,x) =0,\quad (t,x) \in \lbrack 0,T]\times
\partial \Omega ;%
\end{cases}%
\end{equation}
where $b_{1},b_{2}$ are constants and $\Omega \subset R^{n}$ is a bounded domain with
smooth boundary. Let  $E=L^{2}( R)$, $D( A) =H^{2}( \Omega
) \cap  H_{0}^{1}( \Omega ) $ and $A\nu =\Delta \nu $, $%
\forall \nu \in D( A) $. Then (4.4) can be written as
\begin{equation}
\begin{cases}
\widetilde{U_{tt}}( t) +b_{1}A\tilde{U_{t}}( t)
+b_{2}A^{2}\tilde{U}( t) =\widetilde{F( t) }; \\[.3pc]
\tilde{U}( 0) =\Psi _{1},\ \tilde{U_{t}}( 0)
=\Psi _{2}\ .%
\end{cases}
\end{equation}
Through a simple calculus, one may get $\frac{{\rm d}^{k}}{{\rm d}t^{k}}A^{2-k}h(
t) =A^{2-k}\frac{{\rm d}^{k}}{{\rm d}t^{k}}h( t) $, $k=1,2$\ for any $%
h\in C^{2}( [ 0,T]\hbox{\rm :}\ D( A^{2}) ) $. Pazy (p.~211 of
\cite{8}) shows that $A$ is the generator of an analytic
$C_{0}$-semigroup \{$\tilde{T}( t) $\}$_{t\geq 0}$. It is easy to
show that (see e.g., p.~104 of \cite{8}) for any$\ \ g\in L^{2}(
\Omega ) $, the initial boundary value problem
\begin{equation}
\begin{cases}
\dfrac{\partial }{\partial t}( y( t,x) ) =\Delta y( t,x),\quad ( t,x)
\in ( 0,T) \times \Omega;\\[.6pc]
y( 0,x) =g( x),\quad x\in \Omega;\\[.2pc]
y( t,x) =0,\quad ( t,x) [ 0,T] \times \partial \Omega;
\end{cases}
\end{equation}
has a unique solution $y( t,x) =( \tilde{T}( t) g) ( x) $.
However, the initial boundary value problem (4.6) can be solved by
separation of variables method. Its solution can be represented as
\begin{equation*}
y( t,x) =\sum_{k=0}^{\infty }\alpha _{k}\hbox{e}^{\lambda _{k}t}w_{k}( x),
\end{equation*}
where $0>\lambda _{1}\geq \lambda _{2}\geq \lambda _{3}\geq \dots $ are the
eigenvalues of the Laplace equation with Dirichlet boundary conditions and $w_{k}( x)
$ are the corresponding normalized eigenfunctions. Thus
\begin{equation*}
( \tilde{T}( t) g) ( x) =\sum_{k=0}^{\infty }\beta
_{k,g}\hbox{e}^{\lambda _{k}t}w_{k}( x).
\end{equation*}
If the characteristic equation $P( z) =z^{2}+b_{1}z+b_{2}$ of
(4.5) has a root $\widetilde{\alpha _{1}}$ with multiplicity 2,
then $\widetilde{\alpha _{1}}A$\ generates a $C_{0}$-semigroup$\{
\tilde{ T_{1}}( t) \} _{t\geq 0}$ and the solution of (4.5) can be
represented as
\begin{align*}
( \tilde{T}_{1}( t)
g) ( x) =\sum_{k=0}^{\infty }\beta _{k,g}\hbox{e}^{\lambda _{k}%
\widetilde{\alpha _{1}}t}w_{k}( x) \ \hbox{for}\ g\in L^{2}(
\Omega) ,( t,x) \in (0,T) \times \Omega.
\end{align*}
As Example~1, this implies that
\begin{align*}
U( t) =T_{1}( t) \Psi _{1}+tT_{1}( t) ( \phi
_{2}-\widetilde{\alpha _{1}}A\phi _{1}) +\int_{0}^{t}( t-t_{2})
T_{1}( t-t_{2}) \tilde{F}( t_{2})\, \d t_{2}.
\end{align*}
Furthermore, if $\Psi _{1}\in D( A^{3}) $, $\Psi _{2}\in D( A^{2})
$, $\tilde{F}\in C^{2}( [ 0,T] ;L^{2}( \Omega ) ) $, $\tilde{F}(
t) \in D( A^{2}) $ for all $t$ in $( 0,T) $ and $\tilde{F}( 0) \in
D( A) $, then the solution of (4.4) can represented as
\begin{align*}
u( t,x) &=\sum_{k=0}^{\infty }[ \beta
_{k,\Psi _{1}}+t\beta _{k,\Psi _{2}-\widetilde{\alpha _{1}}\Delta \Psi _{1}}%
]\, \hbox{e}^{\lambda _{k}\widetilde{\alpha _{1}}t}w_{k}( x) \\[.4pc]
&\quad\,+\int_{0}^{t}( t-s)
\left[ \sum_{k=0}^{\infty }\beta _{k,f( s,x) }\hbox{e}^{\lambda _{k}%
\widetilde{\alpha _{1}}(t-s)}w_{k}( x) \right]\, \d s.
\end{align*}}
\end{exam}

\section*{Acknowledgement}

The authors would like to thank the referee for useful suggestions
and modifications to this paper.

\end{document}